\documentclass[11pt]{article}
\usepackage{amsfonts,amssymb,color}
\usepackage[mathscr]{eucal}
\usepackage{amsmath, amsthm}
\usepackage{mathrsfs}
\usepackage{amsbsy}
\usepackage{dsfont}
\usepackage{bbm}
\usepackage{wasysym}
\usepackage{stmaryrd}
\usepackage{upgreek}
\usepackage{hyperref}
\input xypic
\xyoption{all}
\begin{document}
\baselineskip = 5mm
\newcommand \primes {{\bf P}}
\newcommand \ZZ {{\mathbb Z}}
\newcommand \FF {{\mathbb F}}
\newcommand \NN {{\mathbb N}}
\newcommand \QQ {{\mathbb Q}}
\newcommand \RR {{\mathbb R}}
\newcommand \CC {{\mathbb C}}
\newcommand \PR {{\mathbb P}}
\newcommand \AF {{\mathbb A}}
\newcommand \VV {{\mathbb V}}
\newcommand \ud {{\Bbbk }}
\newcommand \bgf {{\mathbb K}}
\newcommand \bcA {{\mathscr A}}
\newcommand \bcB {{\mathscr B}}
\newcommand \bcC {{\mathscr C}}
\newcommand \bcD {{\mathscr D}}
\newcommand \bcE {{\mathscr E}}
\newcommand \bcF {{\mathscr F}}
\newcommand \bcG {{\mathscr G}}
\newcommand \bcH {{\mathscr H}}
\newcommand \bcM {{\mathscr M}}
\newcommand \bcN {{\mathscr N}}
\newcommand \bcI {{\mathscr I}}
\newcommand \bcJ {{\mathscr J}}
\newcommand \bcK {{\mathscr K}}
\newcommand \bcL {{\mathscr L}}
\newcommand \bcO {{\mathscr O}}
\newcommand \bcP {{\mathscr P}}
\newcommand \bcQ {{\mathscr Q}}
\newcommand \bcR {{\mathscr R}}
\newcommand \bcS {{\mathscr S}}
\newcommand \bcT {{\mathscr T}}
\newcommand \bcU {{\mathscr U}}
\newcommand \bcV {{\mathscr V}}
\newcommand \bcW {{\mathscr W}}
\newcommand \bcX {{\mathscr X}}
\newcommand \bcY {{\mathscr Y}}
\newcommand \bcZ {{\mathscr Z}}
\newcommand \goa {{\mathfrak a}}
\newcommand \gob {{\mathfrak b}}
\newcommand \goc {{\mathfrak c}}
\newcommand \god {{\mathfrak d}}
\newcommand \gom {{\mathfrak m}}
\newcommand \gon {{\mathfrak n}}
\newcommand \goo {{\mathfrak o}}
\newcommand \gop {{\mathfrak p}}
\newcommand \goq {{\mathfrak q}}
\newcommand \goA {{\mathfrak A}}
\newcommand \goB {{\mathfrak B}}
\newcommand \goT {{\mathfrak T}}
\newcommand \goC {{\mathfrak C}}
\newcommand \goD {{\mathfrak D}}
\newcommand \goM {{\mathfrak M}}
\newcommand \goN {{\mathfrak N}}
\newcommand \goO {{\mathfrak O}}
\newcommand \goP {{\mathfrak P}}
\newcommand \goQ {{\mathfrak Q}}
\newcommand \goS {{\mathfrak S}}
\newcommand \goH {{\mathfrak H}}
\newcommand \ga {{\mathfrak a}}
\newcommand \gb {{\mathfrak b}}
\newcommand \gc {{\mathfrak c}}
\newcommand \gd {{\mathfrak d}}
\newcommand \gm {{\mathfrak m}}
\newcommand \gn {{\mathfrak n}}
\newcommand \gp {{\mathfrak p}}
\newcommand \gq {{\mathfrak q}}
\newcommand \gQ {{\mathfrak Q}}
\newcommand \gP {{\mathfrak P}}
\newcommand \gT {{\mathfrak T}}
\newcommand \gC {{\mathfrak C}}
\newcommand \gD {{\mathfrak D}}
\newcommand \gM {{\mathfrak M}}
\newcommand \gS {{\mathfrak S}}
\newcommand \gH {{\mathfrak H}}
\newcommand \gA {{\rm {A}}}
\newcommand \gB {{\rm {B}}}
\newcommand \catC {{\sf C}}
\newcommand \catD {{\sf D}}
\newcommand \catF {{\sf F}}
\newcommand \catG {{\sf G}}
\newcommand \catE {{\sf E}}
\newcommand \catI {{\sf I}}
\newcommand \catS {{\sf S}}
\newcommand \catW {{\sf W}}
\newcommand \catX {{\sf X}}
\newcommand \catY {{\sf Y}}
\newcommand \catZ {{\sf Z}}
\newcommand \Sets {{\sf Sets}}
\newcommand \Sch {{\sf Sch }}
\newcommand \Funct {{\rm Funct}}
\newcommand \PShv {{\sf PShv}}
\newcommand \Shv {{\sf Shv }}
\newcommand \APSh {{\sf APSh}}
\newcommand \AShv {{\sf AShv }}
\newcommand \Ringspace {{\sf {Ringspace}}}
\newcommand \Nor {{\sf Nor}}
\newcommand \Seminor {{\sf sNor}}
\newcommand \Reg {{\sf Reg}}
\newcommand \Sm {{\sf Sm}}
\newcommand \SmProj {{\sf SmProj}}
\newcommand \NorCon {{\sf NorCon}}
\newcommand \Noe {{\sf Noe}}
\newcommand \LocNoe {{\sf LocNoe}}
\newcommand \Stk {{\sf Stk}}
\newcommand \Mon {{\sf Mon }}
\newcommand \Mod {{\sf Mod}}
\newcommand \Ab {{\sf Ab }}
\newcommand \Ind {{\sf {Ind}}}
\newcommand \Vect {{\sf Vect}}
\newcommand \MM {{\sf MM}}
\newcommand \HS {{\sf HS}}
\newcommand \MHS {{\sf MHS}}
\newcommand \Zar {\rm {Zar}}
\newcommand \Nis {\rm {Nis}}
\newcommand \Nen {\rm {N\acute en}}
\newcommand \cdh {\rm {cdh}}
\newcommand \h {\rm {h}}
\newcommand \fppf {{\rm {fppf}}}
\newcommand \et {\rm {\acute e t}}
\newcommand \CHM {{\sf Chow}}
\newcommand \DM {{\sf DM}}
\newcommand \Ta {{\mathbbm T}}
\newcommand \uno {{\mathbbm 1}}
\newcommand \Le {{\mathbbm L}}
\newcommand \ptr {{\pi _2^{\rm tr}}}
\newcommand \Spec {{\rm {Spec}}}
\newcommand \bSpec {{\bf {Spec}}}
\newcommand \Pic {{\rm {Pic}}}
\newcommand \PicF {{\it {Pic}}}
\newcommand \PicS {{\rm {Pic}}}
\newcommand \Jac {{{J}}}
\newcommand \AlbS {{\rm {Alb}}}
\newcommand \alb {{\rm {alb}}}
\newcommand \NS {{{NS}}}
\newcommand \Corr {{Corr}}
\newcommand \Sym {{\rm {Sym}}} 
\newcommand \Symcl {{\rm {S}}} 
\newcommand \Alt {{\rm {Alt}}}
\newcommand \Prym {{\rm {Prym}}}
\newcommand \HilbF {{\it Hilb}}
\newcommand \HilbS {{\it Hilb}}
\newcommand \Mor {{\rm Mor}}
\newcommand \MorF {{\it Mor}}
\newcommand \MorS {{\it Mor}}
\newcommand \pdiv {{\rm {div}}}
\newcommand \Bl {{\rm {Bl}}}
\newcommand \codim {{\rm {codim}}}
\newcommand \Proj {{\rm {Proj}}}
\newcommand \bProj {{\bf {Proj}}}
\newcommand \Div {{\rm {Div}}}
\newcommand \prim {{\rm {prim}}}
\newcommand \stalk {{\rm st}}
\newcommand \seminorm {{\rm {sn}}}
\newcommand \hc {{\rm hc}} 
\newcommand \BC {{\rm {BC}}}
\newcommand \cnj {{\rm {C}}}
\newcommand \Res {{\rm {Res}}}
\newcommand \Cycl {{\it Cycl }}
\newcommand \PropCycl {{\it PropCycl }}
\newcommand \cycl {{\it cycl }}
\newcommand \PrimeCycl {{\it PrimeCycl }}
\newcommand \PrimePropCycl {{\it PrimePropCycl }}
\newcommand \card {{\rm {card}}}
\newcommand \cone {{\rm {cone}}}
\newcommand \cha {{\rm {char}}}
\newcommand \eff {{\rm {eff}}}
\newcommand \cl {{\rm {cl}}}
\newcommand \tr {{\rm {tr}}}
\newcommand \pr {{\rm {pr}}}
\newcommand \ev {{\rm {ev}}}
\newcommand \interior {{\rm {Int}}}
\newcommand \sep {{\rm {sep}}}
\newcommand \td {{\rm {tdeg}}}
\newcommand \alg {{\rm {alg}}}
\newcommand \im {{\rm im}}
\newcommand \rd {{\rm {red}}}
\newcommand \hl {{\rm h}}
\newcommand \shl {{\rm sh}}
\newcommand \op {{\rm op}}
\newcommand \Hom {{\rm Hom}}
\newcommand \uHom {{\underline {\rm Hom}}}
\newcommand \cHom {{\mathscr H\! }{\it om}}
\newcommand \Ext {{\rm Ext}}
\newcommand \cExt {{\mathscr E\! }{\it xt}}
\newcommand \colim {{{\rm colim}\, }} 
\newcommand \End {{\rm {End}}}
\newcommand \coker {{\rm {coker}}}
\newcommand \id {{\rm {id}}}
\newcommand \van {{\rm {van}}}
\newcommand \spc {{\rm {sp}}}
\newcommand \Ob {{\rm Ob}}
\newcommand \Aut {{\rm Aut}}
\newcommand \cor {{\rm {cor}}}
\newcommand \res {{\rm {res}}}
\newcommand \Gal {{\rm {Gal}}}
\newcommand \PGL {{\rm {PGL}}}
\newcommand \Gr {{\rm {Gr}}}
\newcommand \Tor {{\rm {Tor}}}
\newcommand \Sing {{\rm {Sing}}}
\newcommand \spn {{\rm {span}}}
\newcommand \univ {{\rm {\, univ}}}
\newcommand \Nm {{\rm {Nm}}}
\newcommand \fin {{\rm {f}}}
\newcommand \inv {{\rm {inv}}}
\newcommand \even {{\rm {even}}}
\newcommand \md {{\rm {mod}\, }}
\newcommand \sg {{\Sigma }}
\newcommand \ind {{\rm {ind}}}
\newcommand \Gm {{{\mathbb G}_{\rm m}}}
\newcommand \trdeg {{\rm {tr.deg}}}
\newcommand \con {\rm {conn}}
\newcommand \sv {{\rm {sv}}}
\newcommand \sing {{\rm {sing}}}
\newcommand \tame {\rm {tame }}
\newcommand \eq {{\rm {eq}}}
\newcommand \length {{\rm {length}}}
\newcommand \ord {{\rm {ord}}}
\newcommand \shf {{\rm {a}}}
\newcommand \spd {{\rm {s}}}
\newcommand \glue {{\rm {g}}}
\newcommand \equi {{\rm {equi}}}
\newcommand \ab {{\rm {ab}}}
\newcommand \add {{\rm {ad}}}
\newcommand \Fix {{\rm {Fix}}}
\newcommand \pty {{\mathbf P}}
\newcommand \type {{\mathbf T}}
\newcommand \trp {{\rm {t}}}
\newcommand \cat {{\rm {cat}}}
\newcommand \deop {{\Delta \! }^{op}\, }
\newcommand \defect {{\rm {def}}}
\newcommand \aff {{\rm {aff}}}
\newcommand \Const {{\rm {Const}}}
\newcommand \num {{\rm {num}}}
\newcommand \conv {{\it {cv}}}
\newcommand \nil {{\rm {nil}}}
\newcommand \rat {{\rm rat}}
\newcommand \gen {{\rm gen}}
\newcommand \SG {{\rm SG}}
\newcommand \tors {{\rm {tors}}}
\newcommand \coeq {{{\rm coeq}\, }}
\newcommand \supp {{\rm Supp}}
\newcommand \sm {{\rm sm}}
\newcommand \reg {{\rm reg}}
\newcommand \nor {{\rm nor}}
\newcommand \noe {{\rm Noe}}
\newcommand \var {{\rm var}}
\newcommand \norm {{\rm {nor}}}
\newcommand \tp {{\rm {tp}}}
\newcommand \st {{\rm {st}}}
\newcommand \Gys {{\rm {Gys}}}
\newcommand \Fr {{\rm {Fr}}}
\newcommand \hol {{\rm {h}}}
\newcommand \bideg {{\rm {bideg}}}
\newcommand \loc {{\rm {loc}}}
\newcommand \Rat {{\mathscr Rat}}
\newcommand \tdf {{\mathbf {t}}}
\newcommand \term {{*}}
\newcommand \znak {{\natural }}
\newcommand \znakk {{\sharp }}
\newcommand \znakkk {{\flat }}
\newcommand \qand {{\quad \hbox{and}\quad }}
\newcommand \qqand {{\quad \quad \hbox{and}\quad \quad }}
\newcommand \qqqand {{\quad \quad \quad \hbox{and}\quad \quad \quad }}
\newcommand \heither {{\hbox{either}\quad }}
\newcommand \qor {{\quad \hbox{or}\quad }}
\newcommand \qqor {{\quad \quad \hbox{or}\quad \quad }}
\newcommand \lra {\longrightarrow}
\newcommand \hra {\hookrightarrow}
\def\blue {\color{blue}}
\def\red {\color{red}}
\def\green {\color{green}}
\newtheorem{theorem}[subsubsection]{Theorem}
\newtheorem{lemma}[subsubsection]{Lemma}
\newtheorem{corollary}[subsubsection]{Corollary}
\newtheorem{proposition}[subsubsection]{Proposition}
\newtheorem{remark}[subsubsection]{Remark}
\newtheorem{definition}[subsubsection]{Definition}
\newtheorem{conjecture}[subsubsection]{Conjecture}
\newtheorem{example}[subsubsection]{Example}
\newtheorem{question}[subsubsection]{Question}
\newtheorem{comment}[subsubsection]{Comment}
\newtheorem{assumption}[subsubsection]{Assumption}
\newtheorem{fact}[subsubsection]{Fact}
\newtheorem{crucialquestion}[subsubsection]{Crucial Question}
\newtheorem{claim}[subsubsection]{Claim}
\newtheorem{terminology}[subsubsection]{Terminology}
\newtheorem{aim}[subsubsection]{Aim}
\newenvironment{pf}{\par\noindent{\em Proof}.}{\hfill\framebox(6,6)
\par\medskip}
\title{\bf The Bloch conjecture}


\author{\sc Vladimir Guletski\u \i }


\maketitle

\begin{abstract}
\noindent We prove the conjecture stated by Spencer Bloch in 1975 and saying that the Albanese kernel of a smooth projective surface is $0$, provided its second cohomology group is algebraic.
\end{abstract}
















\tableofcontents

\section{Introduction}
\label{intro}

Let $k$ be an algebraically closed field, and let $X$ be a smooth projective variety over $k$. The goal of this manuscript is to discuss the Chow group $CH_0(X)$ of $0$-cycle modulo rational equivalence on $X$, which is an important invariant of $X$. 

Let $A_0(X)$ be the subgroup in $CH_0(X)$ generated by $0$-cycles of degree $0$, and let $T(X)$ be the kernel of the Albanese homomorphism from $A_0(X)$ onto the group of closed points of the Albanese variety $\AlbS _{X/k}$ of the variety $X$. 

If $X$ is a smooth projective curve over $k$, then $T(X)=0$, and the group $A_0(X)$ is isomorphic to the Jacobian variety of the curve $X$. When $k$ is $\CC $, this is the content of the classical theorem due to Abel and Jacobi proven in the first half of $19$th century. 

When $X$ is a smooth projective surface over $k$, the behaviour of rational equivalence of $0$-cycles on $X$ may be very different from what could be expected based on the intuition taken from curves. In the turning paper \cite{Mumford} Mumford proved that when $k=\CC $ and the geometric genus $p_g(X)>0$, then $T(X)$ is huge, and $A_0(X)$ cannot be parametrized by an abelian variety, in any reasonable sense. 

On the other hand, based on some Hodge and $K$-theoretic intuition, Spencer Bloch has conjectured that $T(X)$ is $0$ for a smooth projective surface $X$ over $\CC $, provided $p_g(X)=0$, see \cite{BCinitial} and \cite{BlochLectures}. Notice that, since $k=\CC $, due to Hodge theory and Lefschetz $(1,1)$-theorem, the condition $p_g=0$ is equivalent to saying that the second transcendental cohomology group $H^2_{\tr }(X,\QQ )$ vanishes. 

The status of BC is as follows. If the Kodaira dimension of $X$ is strictly smaller than $2$, it was proven in \cite{BKL} by ad hoc methods. If $X$ is of general type over $\CC $, the condition $H^2_{\tr }(X,\QQ )=0$ implies that $X$ is regular, i.e. $H^1(X,\QQ )=H^3(X,\QQ )=0$, and then BC simply asserts that any two points on $X$ are rationally equivalent to each other. This is the hard case of BC, unknown for the most part, except for several cases: the surfaces studied by Barlow in \cite{Barlow1} and \cite{Barlow2}, the classical Godeaux surface is considered in \cite{InoseMizukami}, and general Godeaux case was studied by Voisin in \cite{Sur les 0-cycles} and \cite{VoisinVariations}, the Catanese surfaces and general Barlow surfaces solved by Voisin in \cite{VoisinCataneseBarlowSurfaces}, and finite quotients of products of curves, proven by Kimura in \cite{Kimura}. 


Thinking in $l$-adic terms, where $l$ is a prime different from $\cha (k)$, the condition $H^2_{\tr }(X,\QQ _l)=0$ should be strengthened by imposing the condition $H^1(X,\QQ _l)=H^3(X,\QQ _l)=0$, which is automatic in $0$ characteristic, but, surprisingly, in positive characteristic there may be surfaces of general type with $H^2_{\tr }=0$ and non-trivial $H^1$ and $H^3$, see \cite{LiedtkeGodeaux}. This extended condition can be restated by saying that the total cohomology group $H^*(X,\QQ _l)$ is algebraic. 

The goal of this note is to prove the following

\begin{itemize}
\item[]{}
{\rm THEOREM A.} {\it
Let $k$ be an algebraically closed field of characteristic $0$, and let $X$ be a smooth projective surface over $k$. If $H^*(X,\QQ _l)$ is algebraic, then $T(X)=0$.
}
\end{itemize}

The same is true in positive characteristic, when $X$ has a model over the algebraic closure $\bar \FF _p$ of the finite field $\FF _p=\ZZ /p$.

\begin{itemize}
\item[]{}
{\rm THEOREM B.} {\it
Let $k$ be an algebraically closed field of characteristic $p>0$, $p\neq 2$, and let $X$ be a smooth projective surface over $k$. Assume that $X$ has a model over $\bar \FF _p$. If $H^*(X,\QQ _l)$ is algebraic, then $T(X)=0$.
}
\end{itemize}

\medskip

{\sc Acknowledgements.} I am grateful to Kalyan Banerjee, Kestutis \v Cesnavi\v cius, Ivan Cheltsov, Jean-Louis Colliot-Th\'el\`ene, Sergei Gorchinskiy, Joe Palacios Baldeon, Aleksandr Pukhlikov, Will Sawin, Alexander Tikhomirov, Claire Voisin and Bo Zhang for useful discussions on the internet and offline. 










\bigskip

\section{The conjecture}
\label{description&philosophy}

\subsection{Classical formulation}
\label{nonmotivic}

Let $k$ be an algebraically closed field, and let $X$ be a variety over $k$. Consider the Chow group $CH_0(X)$ of $0$-cycles modulo rational equivalence on $X$, and let $A_0(X)$ be its subgroup generated by $0$-cycles of degree $0$ on $X$. 

Assume that $X$ is smooth and projective, and let $\AlbS _{X/k}$ be the Albanese variety of the variety $X$, i.e. the dual to the connected component $\PicS ^0_{X/k}$ of the Picard scheme $\PicS _{X/k}$ of $X$ over $k$. Consider the Albanese homomorphism
  $$
  \alb _X:A_0(X)\to \AlbS _{X/k}(k)
  $$
associated to the variety $X$. Its kernel 
  $$
  T(X)=\ker (\alb _X)
  $$ 
is an important object of study in algebraic geometry. 

The homomorphism $\alb _X$ is an isomorphism on $m$-torsion subgroups, where $m$ is coprime to $\cha (k)$. This is Roitman's theorem, see \cite{Roitman} and \cite{MilneRoitmanThm} 

Let $\Omega ^q_{X/k}$ be the sheaf of K\"ahler differentials of degree $q$ on $X$, and let 
  $$
  h^{p,q}(X)=\dim _k(H^p(X,\Omega ^q_{X/k}))\; .
  $$
be the $(p,q)$-Hodge number of the variety $X$. If $n$ is the dimension of $X$, the number 
  $$
  p_g(X)=h^{0,n}(X)
  $$ 
is called the geometric genus of the variety $X$. 

The original Bloch's conjecture (BC) can be stated as follows, see page 408 in \cite{BCinitial} and page 0.5 in \cite{BlochLectures}. 

\begin{conjecture}
\label{BCclassical}
Let $X$ be a connected smooth projective surface over $\CC $. Then $p_g(X)=0$ implies $T(X)=0$.
\end{conjecture}

As we mentioned in Introduction, the status of BC is that it is proven for all surfaces with Kodaira dimension $<2$, and for a few surfaces of general type. Therefore, further below we only focus on BC for surfaces of general type. Special type will be systematically ignored. 

We would like to generalize Conjecture \ref{BCclassical} to fields of arbitrary characteristic. Let $k$ be an arbitrary algebraically closed field, and fix a prime number $l$ different from the characteristic of $k$. For any algebraic scheme $X$ over $k$, let
  $$
  H^i(X,\ZZ _l)=
  \lim _{\nu }H^i(X,\ZZ /l^{\nu }\ZZ )
  $$
be the $i$-th \'etale $l$-adic cohomology group of the scheme $X$ with coefficients in $\ZZ _l$, and let  
  $$
  H^i(X,\QQ _l)=
  H^i(X,\ZZ _l)\otimes _{\ZZ _l}\QQ _l
  $$
be the corresponding vector space over $\QQ _l$. After a non-canonical choice of an isomorphism 
  $$
  \ZZ _l(1)\simeq \ZZ _l\; ,
  $$
the contravariant functor
  $$
  H^*:X\mapsto H^*(X,\QQ _l)=\oplus _iH^i(X,\QQ _l)
  $$
gives us a Weil cohomology theory over $k$, see \cite{Kleiman}. 

Assume $X$ is equidimensional of dimension $n$ over $k$, and let 
   $$
   CH^i(X)=CH_{n-i}(X)
   $$
be the Chow group of codimension $i$ (respectively, dimension $n-i$) algebraic cycles modulo rational equivalence on $X$. Write also $CH^i(X,\QQ )$ for the Chow group with coefficients in $\QQ $, and similarly for $\QQ _l$.

When $X$ is smooth over $k$, for any integer $0\leq i\leq n$, we have the cycle class homomorphism
  $$
  \cl _X^i:CH^i(X,\QQ _l)\to H^{2i}(X,\QQ _l)
  $$
from the Chow group with coefficients in $\QQ _l$ to $l$-adic cohomology, where the Tate twist is ignored (see $\S 9$, Ch. VI in \cite{MilneEC} or \cite{SGA4.5}, Expos\'e 4).

Let now $X$ be a smooth projective surface over $k$, and let
  $$
  H^2_{\alg }(X,\QQ _l)=\im (\cl ^1_X)
  $$
be the image of the second cycle class homomorphism sending divisors to their cohomology classes on $X$, i.e. the algebraic part in $H^2(X)$, and let 
  $$
  H^2_{\tr }(X,\QQ _l)=H^2(X,\QQ _l)/H^2_{\alg }(X,\QQ _l)
  $$
be the quotient group of $H^2$ by its algebraic part, i.e. the transcendental part of $H^2(X)$, see p. 1.22 in \cite{BlochLectures}. We will say that $H^*_{\et }(X,\QQ _l)$ is algebraic if 
  $$
  H^2_{\tr }(X,\QQ _l)=0
  $$ 
and 
  $$
  H^1(X,\QQ _l)=H^3(X,\QQ _l)=0\; .
  $$

Notice that, if 
  $$
  \cha (k)=0\; ,
  $$ 
then
  $$
  p_g(X)=0\; \; \Leftrightarrow \; \; H^2_{\tr }(X,\QQ _l)=0
  $$
due to Hodge theory. Moreover, we have the following

\begin{lemma}
\label{ehhh}
Let $k$ be an algebraically closed field, let $X$ be a smooth projective surface over $k$, and assume that $X$ is of general type. If $p_g(X)=0$, then $H^1(X,\QQ _l)=H^3(X,\QQ _l)=0$. 
\end{lemma}

\begin{pf}
In arbitrary characteristic, $\dim (\AlbS _{X/k})\leq h^{0,1}(X)$, see again page 240 in \cite{Liedtke}, and if the surface $X$ is minimal and of general type, then $p_g(X)=0$ implies $h^{0,1}(X)=0$ by Proposition 8.9 on page 264 in loc.cit, whence $\AlbS _{X/k}=0$. It follows that $H^1(X,\QQ _l)=0$, and by the Poicar\'e duality, $H^3(X,\QQ _l)=0$.
\end{pf}

However, if $\cha (k)=p>0$ the connection between the geometric genus and transcendental cohomology is more subtle, see \cite{LiedtkeGodeaux}. Since BC is proven for surfaces of Kodaira dimension $<2$, see \cite{BKL}, the best assumption for BC is that the total cohomology group $H^*(X,\QQ _l)$ is algebraic. This is equivalent to saying that $H^1=H^3=0$ and $H^2_{\tr }=0$ (see Appendix to Lecture 2 in \cite{BlochLectures}). 

\begin{conjecture}
\label{BCgeneralized}
Let $k$ be an algebraically closed field, and let $X$ be a smooth projective surface over $k$. If $H^*(X,\QQ _l)$ is algebraic, then $T(X)=0$.
\end{conjecture}

Notice that if $X$ is of general type and $H^*(X,\QQ _l)$ is algebraic, vanishing of $T(X)$ is simply equivalent to saying that any two closed points are rationally equivalent to each other on the surface $X$. 

\subsection{Finite-dimensional motives}
\label{motivic-I}

Let $k$ be an algebraically closed field, and let $\CHM (k,\QQ )$ be the category of Chow motives over $k$ with coefficients in $\QQ $. Morphisms in $\CHM (k,\QQ )$ are correspondences, i.e. classes in $CH^*(X\times Y,\QQ )$, for smooth projective $X$ and $Y$ over $k$. 

The details about this category can be found in \cite{Scholl}. All we need at the moment is that it is pseudo-abelian, tensor, rigid and $\QQ $-linear in the sense of \cite{DeligneMilne}, and the tensor product is induced by the product of smooth projective varieties over $k$. 

In any such a category one can work with wedge and symmetric powers of objects and morphisms, and hence use the machinery of finite-dimensional objects, as it was revolutionary introduced in \cite{Kimura} (see also a vast development in \cite{AK}). 

An object $M$ in $\CHM (k,\QQ )$ is said to be evenly finite-dimensional if there exists a positive integer $n$, such that 
  $$
  \wedge ^nM=0
  $$ 
in $\CHM (k,\QQ )$. The object $M$ is oddly finite-dimensional if 
  $$
  \Sym ^nM=0
  $$ 
in $\CHM (k,\QQ )$, for some $n$. An object $M$ in $\CHM (k,\QQ )$ is finite-dimensional, if it can be decomposed in to a direct sum 
  $$
  M=M_+\oplus M_-\; ,
  $$
where $M_+$ is evenly and $M_-$ is oddly finite-dimensional.

Any Weil cohomology theory $H^*$ can be extended to a functor on $\CHM (k,\QQ )$, see \cite{Scholl}. If $M$ is a Chow motive over $k$, then $H^*(M)$ decomposes in to the even 
  $$
  H^+(M)=\oplus _{i\in \ZZ }H^{2i}(M)
  $$ 
and odd 
  $$
  H^-(M)=\oplus _{i\in \ZZ }H^{2i+1}(M)
  $$ 
parts, since the Weil cohomology theory is $\ZZ $-graded. If $M$ is finite-dimensional, then
  $$
  H^*(M_+)=H^+(M)\; ,
  $$
  $$
  H^*(M_-)=H^-(M)\; ,
  $$
and, clearly, 
  $$
  \wedge ^{b_++1}M_+=0\; ,
  $$
  $$
  \Sym ^{b_-+1}M_-=0\; ,
  $$
where $b_+$ and $b_-$ are the even and odd Betti numbers of the motive $M$.

All these results can be found in \cite{Kimura}.

In particular, if $X$ is a smooth projective surface with $b_1(X)=0$ over $k$, then $b_3(X)=0$ by Poincar\'e duality in the Weil cohomology theory $H^*$, and then finite-dimensionality of the motive $M(X)$ of the surface $X$ is equivalent to even finite-dimensionality of $M(X)$, and in turn the latter is equivalent to vanishing 
  $$
  \wedge ^{b_+(X)+1}M(X)=0\; .
  $$

If $L$ is a field extension over $k$, then we have a natural tensor functor from $\CHM (k,\QQ )$ to $\CHM (L,\QQ )$, sending a motive $M$ to its scalar extension $M_L$. 

\begin{lemma}
\label{scalars}
If $M$ is finite-dimensional, then $M_L$ is finite-dimensional. If $M_L$ is evenly (oddly) finite-dimensional, then so is $M$. If $M$ is the motive of a smooth projective surface over $k$, $M$ is finite-dimensional if and only if $M_L$ is finite-dimensional. 
\end{lemma}

\begin{pf}
The scalar extension functor from $\CHM (k,\QQ )$ to $\CHM (L,\QQ )$ is obviously additive and tensor. Hence, it preserves finite-dimensional objects. Since we look at motives with coefficients in $\QQ $, if $M_L$ is evenly (oddly) finite-dimensional, then $M$ is evenly (oddly) finite-dimensional due to Lemma 3 on page 1.21 in \cite{BlochLectures}. If $M$ is the motive of a smooth projective surface over $k$, we can split out the Picard and Albanese motives from $M$, which correspond to only odd cohomology, and so decompose the diagonal in to even and odd parts, see \cite{Murre} or \cite{Scholl}. The needed equivalence then follows from the equivalence in the even and odd cases.
\end{pf}

The meaning of Lemma \ref{scalars} is that finite-dimensionality of the motive of a surface does not depend on the field of definition, and any hypothesis, which can be restated in terms of motivic finite-dimensionality, is ground field independent too. 

Let $\CHM (k,\QQ )^{\ab }$ be the full tensor subcategory in $\CHM (k,\QQ )$ generated (additively and tensorly) by motives of curves. In other words, $\CHM (k,\QQ )^{\ab }$ is the full tensor subcategory of motives of abelian type in $\CHM (k,\QQ )$. 

\begin{proposition}
\label{Kimuraab}
For any algebraically closed field $k$, all motives in $\CHM (k,\QQ )^{\ab }$ are finite-dimensional. 
\end{proposition}

\begin{pf}
See Corollary 4.4 on page 184 in \cite{Kimura}, and use the standard categorical properties of finite-dimensional motives.
\end{pf}

Morphisms of Chow motives are classes of algebraic cycles modulo rational equivalence on varieties over a field. As such, they can be equivalent to $0$ with regard to any adequate equivalence relation on algebraic cycles. 

The key result in the theory of finite-dimensional motives is due to S.-I. Kimura:

\begin{theorem}
\label{Kimurakey1}
Let $k$ be an algebraically closed field, and let $M$ be a finite-dimensional object of the category $\CHM (k,\QQ )$. If $f:M\to M$ is an endomorphism of $M$ and $f$ is numerically equivalent to $0$, then $f$ is nilpotent in the associative algebra $\End (M)$. 
\end{theorem}

\begin{pf}
See Proposition 7.5 on pp 192 - 193 in \cite{Kimura}.
\end{pf}

As a consequence of this theorem, Kimura also proved

\begin{corollary}
\label{Kimurakey2}
Let $k$ be an algebraically closed field, and let $X$ be a smooth projective variety over $k$, such that its Chow motive $M(X)$ is finite-dimensional. If the total cycle class homomorphism
  $$
  \cl _X^*:CH^*(X,\QQ _l)\to H^*(X,\QQ _l)
  $$
is surjective, i.e. $H^*(X,\QQ _l)$ is algebraic, then it is also an isomorphism of $\QQ _l$-vector spaces. 
\end{corollary}

\begin{pf}
See Proposition 7.6 on p 193 in \cite{Kimura}. 
\end{pf}

It means, in particular, that if $M(X)$ is finite-dimensional and $H^*(X,\QQ _l)$ is algebraic, for a smooth projective surface $X$ over $k$, then $T(X)=0$, i.e. Bloch's conjecture holds for $X$.

The converse is also true: if BC holds for $X$, then $M(X)$ is finite-dimensional, see \cite{GP2} and apply Lemma \ref{scalars}. Thus, Bloch's conjecture, for $X$ with $H^*(X,\QQ _l)$ algebraic, is actually equivalent to finite-dimensionality of the motive $M(X)$.

\subsection{Transcendence degree of $0$-cycles}
\label{motivic-II}

Let $k$ be an algebraically closed field, let $X$ be an equi-dimensional variety of dimension $n$ over $k$, and let $Z$ be an integral closed subscheme of codimension $n$ on $X\times X$, i.e. a prime correspondence of degree $0$ from $X$ to $X$. Following \cite{Barbieri-Viale}, we say that $Z$ is  balanced on the left (respectively, on the right) if there exists an equi-dimensional Zariski closed subscheme $Y\subset X$ with $\dim (Y)<n$, such that $Z$ is contained in $Y\times X$ (respectively, in $X\times Y$). A codimension $n$ cycle $Z$ is balanced on $X\times X$ if it is a sum of prime cycles balanced on the left or right. A class in $CH^n(X\times X)$ is balanced if it can be represented by a balanced cycle on $X\times X$, and the same with coefficients in $\QQ $. 

It is essential for what follows that the subspace generated by balanced cycles is a two-sided ideal in the associative algebra 
  $$
  CH^n(X\times X,\QQ )
  $$
of degree $0$ correspondences with coefficients in $\QQ $ on $X$, see p. 309 in \cite{Fulton}.

Let $k\subset K$ be a field extension. A $K$-point on $X$ is a morphism of schemes 
  $$
  P:\Spec (K)\to X
  $$
over $\Spec (k)$. Let $\zeta _P$ be the image of the unique point in $\Spec (K)$ with respect to the morphism $P$, and let $\kappa (\zeta _P)$ be the residue field of the point $\zeta _P$ on the scheme $X$. By definition, the transcendence degree of the point $P$ over $k$ is the transcendence degree of the field $\kappa (\zeta _P)$ over $k$:
  $$
  \trdeg (P/k)=\trdeg (\kappa (\zeta _P)/k)\; .
  $$
The morphism $P$ induces the field embedding 
  $$
  \kappa (\zeta _P)\subset K\; ,
  $$
whence
  $$
  \trdeg (P/k)\leq \trdeg(K/k)\; .
  $$

Let 
  $$
  Z=\overline {\{ \zeta _P\} }
  $$
be the Zariski closure of the schematic point $\zeta _P$ in $X$. Then $Z$ is a closed integral subscheme in $X$, and
  $$
  \trdeg (P/k)=\dim(Z)\; .
  $$
Therefore, $\trdeg (P/k)$ is the minimal possible dimension of a subvariety $Z$ in $X$ over $k$, such that the point $P$ is in $Z_K$. In particular, that
  $$
  \trdeg (P/k)\leq \dim(X)\; .
  $$

Let $\ud $ be a universal domain over the field $k$, i.e. an algebraically closed field extension of infinite transcendence degree over $k$, and let
  $$
  X_{\ud }=X\times _k\ud 
  $$
be the scalar extension for $X$. The transcendence degree
  $$
  \trdeg (\alpha /k)
  $$ 
of a $0$-cycle class 
  $$
  \alpha \in CH^n(X_{\ud })
  $$ 
over $k$ is the minimal natural number $t$, such that there exists a $0$-cycle
  $$
  Z=\sum _im_iP_i
  $$
on $X_{\ud }$, representing the class $\alpha $, such that 
  $$
  \trdeg (P_i/k)\leq t
  $$
for all $i$.

\medskip

The following theorem is a development of the result in \cite{BlochSrinivas}.

\begin{theorem}
\label{maintrdeg}
Let $X$ be an irreducible smooth projective variety of dimension $n$ over $k$. The following conditions are equivalent:

\begin{enumerate}

\item[]{(i)}
the class of the diagonal $\Delta _X$ is balanced in $CH^n(X\times X,\QQ )$;

\item[]{(ii)}
the Chow motive $M(X)$ is isomorphic to a direct summand of the motive
  $$
  M(Y_1)\oplus (M(Y_2)\otimes \Le ^{n-e})\; ,
  $$
where $\Le ^{n-e}$ is the tensor power of the Lefschetz motive, $Y_1$ and $Y_2$ are equidimensional smooth projective varieties over $k$ whose dimensions are strictly less than $n$, and $e$ is the dimension of the variety $Y_2$;

\item[]{(iii)}
any element $\alpha \in CH^n(X_{\ud },\QQ )$ satisfies
  $$
  \trdeg (\alpha /k)<n\; ;
  $$

\item[]{(iv)}
there exists a closed point $P\in X_{\ud }$, such that 
  $$
  \trdeg (P/k)=n
  $$ 
and
  $$
  \trdeg ([P])<n\; ,
  $$
where $[P]$ is the class of the point $P$ in $CH^d(X_{\ud },\QQ )$.

\end{enumerate}

\end{theorem}

\begin{pf}
This is Theorem 3.6 on page 564 in \cite{trdegzerocycles}.
\end{pf}

As motives of curves are finite-dimensional, Proposition \ref{Kimuraab} and Corollary \ref{Kimurakey2} guarantee that to prove BC, for $X$ with algebraic $H^*(X,\QQ _l)$, it suffices to prove any of the four equivalent items in Theorem \ref{maintrdeg}.

\subsection{An outline of the approach}
\label{outline}


Let $X$ be a smooth projective surface over $k$. Fix a closed immersion of $X$ into some projective space, and use Segre to embed $X\times X$ in to the bigger projective space $\PR $ over $k$. Let $V$ be the section of $X\times X$ by a general hypersurface of sufficiently big degree containing the diagonal $\Delta _X$ in $\PR $. The singular locus $V_{\sing }$ is a finite number of isolated ODPs on $V$ by Proposition 1.4 on page 11 in Bloch's thesis \cite{BlochThesis}. If $k=\CC $, by the result of Brevik and Nollet, see Theorem 1.1 on page 108 in \cite{BrevikNollet}, 
       $$
       CH^1(V)=\langle \Delta _X\rangle \oplus I\; ,
       $$
where $I$ is the image of the restriction homomorphism 
  $$
  CH^1(X\times X)\to CH^1(V)\; .
  $$ 


We shall find a surface $F$ in $V$, such that 
  $$
  F\cap V_{\sing }=\emptyset \; ,
  $$
and $F$ projects dominantly on both factors of $X\times X$. The class of $F$ is then a linear combination of the class of $\Delta _X$ and the restrictions of divisors from $X\times X$ to $V$, by Brevik-Nollet's theorem. If the irregularity of $X$ is $0$, 
  $$
  CH^1(X\times X)=CH^1(X)\oplus CH^1(X)\; ,
  $$
and all such restrictions are balanced. 

If $\Delta _X$ contributes trivially in to the class of $F$ in $CH^1(V)$, the surface $F$ is balanced modulo rational equivalence on $X\times X$. Since $F$ projects dominantly on both factors, it gives rise to a point $P$ of transcendence degree $2$ over $k$ on $X_{\ud }$. As balancing of $F$ is equivalent to saying that 
   $$
   \trdeg ([P]/k)<2\; ,
   $$
 the diagonal $\Delta _X$ is balanced by Theorem \ref{maintrdeg}, and we are done.

The second case is when the contribution of $\Delta _X$ in to the class of $F$ in $CH^1(V)$ is not trivial. Then $\Delta _X$ is rationally equivalent on $V$ to a linear combination of $F$ and the restrictions of divisors from $X\times X$ to $V$, which are balanced provided $X$ is regular. Since $F$ does not meet $V_{\sing }$, we get an honest cohomological class of $F$ in $H^2(V,\QQ _l)$ by excision. This class is balanced modulo homological equivalence by Weak Lefschetz, provided 
  $$
  H^2_{\tr }(X,\QQ _l)=0\; .
  $$

Resolving ODPs on $V$ and applying Matsusaka's theorem, we obtain that $F$ is balanced modulo algebraic equivalence on $V$, and so on $X\times X$. The Voevodsky-Voisin nilpotency theorem yields that $\Delta _X$ minus a balanced correspondence is nilpotent in the associative algebra $CH^0(X,X,\QQ )$. 

The binomial expansion completes the proof.

\section{Proof of Theorem A}
\label{proofA}

\subsection{Nilpotency theorems}

Let $\bcC $ be a rigid tensor category in the sense of \cite{DeligneMilne}. For any morphism 
  $$
  f:M\to N
  $$ 
in $\bcC $, consider its $n$-th tensor power
  $$
  f^{\otimes n}:M^{\otimes n}\to N^{\otimes n}
  $$
as an element in the vector space 
  $$
  \Hom _{\bcC }(M^{\otimes n},N^{\otimes n})\; .
  $$ 
If
  $$
  f^{\otimes n}=0\; ,
  $$
we say that $f$ is smash-nilpotent in the category $\bcC $.

And, as usual, if $M=N$, we say that $f:M\to N$ is nilpotent, if 
  $$
  f^n=0\; ,
  $$ 
where the power is self-composing of $f$ in $\bcC $.

\begin{lemma}
\label{nilponilpo}
In any rigid tensor category, a smash-nilpotent endomorphism is nilpotent.
\end{lemma}

\begin{pf}
See Lemma 2.7 on page 191 in \cite{VoevodskyNilp}, and Lemma 7.4.2 on page 162 in \cite{AK}. 
\end{pf}

Now, let $k$ be an algebraically closed field, and let $X$ be a smooth projective variety over $k$. If $A$ is an algebraic cycle of dimension $d$ on $X$, for any positive integer $n$ we can consider its exterior product $A^{\times n}$ as an algebraic cycle of dimension $dn$ on $X^{\times n}$, and the same for cycle classes in Chow groups, see pages 24 - 25 in \cite{Fulton}. A cycle class 
  $$
  \alpha \in CH_d(X,\QQ )
  $$ 
is said to be smash-nilpotent if 
  $$
  \alpha ^{\times n}=0
  $$ 
in $CH_{dn}(X^{\times n},\QQ )$. 

Looking at elements of Chow groups with coefficients in $\QQ $ as morphisms in $\CHM (k,\QQ )$, we see that exterior products of cycle classes in Chow groups are particular cases of tensor products of morphisms in tensor rigid categories, and smash-nilpotency of cycle classes is smash-nilpotency in $\CHM (k,\QQ )$.

In \cite{VoevodskyNilp} Voevodsky prove the following important 

\begin{theorem}
\label{VNT1}
Let $\alpha $ be a cycle class in $CH_d(X,\QQ )$. If $\alpha $ is algebraically equivalent to $0$, it is smash-nilpotent in the sense above. 
\end{theorem}

\begin{pf}
See Corollary 3.2 on page 192 in \cite{VoevodskyNilp}. 
\end{pf}

Theorem \ref{VNT1} and Lemma \ref{nilponilpo} yield the following

\begin{corollary}
\label{VNT2}
Let $f:M\to M$ be an endomorphism in $\CHM (k,\QQ )$. If $f$ is algebraically equivalent to $0$, then it is nilpotent in $\End (M)$.
\end{corollary}

\begin{pf}
See Corollary 3.3 on page 193 in \cite{VoevodskyNilp}. 
\end{pf}

Similar results were independently obtained by C. Voisin in \cite{VoiN}.

 
\subsection{Bloch's Bertini type theorem}
\label{BBT}

Let $k$ be an algebraically closed field, and let $X$ be a smooth projective surface over $k$. Let also $V$ be an irreducible projective threefold in $X\times X$, which is dominant with regard to both projections from $X\times X$ on to $X$. For any integral closed subscheme $Z$ of dimension $2$ in $V$, we will say that $Z$ is balanced on $V$ if $Z$ is balanced on $X\times X$ in the sense of Section \ref{motivic-II}. If $Z$ is dominant on both factors, we will say sometimes that $Z$ is fat in $X\times X$. An algebraic cycle $Z$ of dimension $2$ on $V$ is balanced on $V$, if all its prime components are balanced on $V$.


For short, let
  $$
  W=X\times X\; ,
  $$
fix a closed immersion of $X$ in to some projective space, and use the Segre embedding to embed $W$ in to a bigger projective space $\PR $ over the ground field $k$. Then the diagonal
  $$
  \Delta _X\subset W\subset \PR 
  $$ 
is also in $\PR $. 

Let $\bcI $ be the sheaf of ideals in $\bcO $ defining the closed subscheme $\Delta _X$ in $\PR $. For a general hypersurface $H$ in the linear system 
  $$
  H^0(\PR ,\bcI (d))\; ,
  $$
when the degree $d$ is sufficiently big, the intersection 
  $$
  V=W\cap H
  $$
is irreducible, containing $\Delta _X$, and the singular locus $V_{\sing }$ of the variety $V$ is a finite collection 
  $$
  \{ P_1,\ldots ,P_n\} 
  $$
of isolated ordinary double points, all supported on $\Delta _X$. This is Proposition 1.2 on pages 9 - 10 and Proposition 1.4 on page 11 in Bloch's thesis \cite{BlochThesis}.
  
\subsection{A fat surface away of ODPs}
\label{fat}

Choose a fix a Noether normalization morphism
  $$
  X\to \PR ^2
  $$
over $k$, and let
  $$
  \nu :X\times X\to \PR ^2\times \PR ^2
  $$
be the fibred square of it over $k$. The scheme-theoretic image
   $$
   T_0=\nu (V)
   $$
is an irreducible threefold in $\PR ^2\times \PR ^2$ containing the diagonal $\Delta _{\PR ^2}$ and, therefore, dominant on both factors in $\PR ^2\times \PR ^2$. 

Let $T_1$ be a general linear section of the Segre fourfold $\PR ^2\times \PR ^2$ in $\PR ^8$. Then $T_1$ is smooth irreducible by Bertini's theorem, and as $T_1$ is a $(1,1)$-divisor in $\PR ^2\times \PR ^2$, it is dominant on both factors too. 

Now, since $T_0$ and $T_1$ are dominant on both factors, their set-theoretical intersection
  $$
  S=T_0\cap T_1
  $$ 
is dominant on both factors, because any two plane curves have non-empty intersection in $\PR ^2$. Moreover, if $T_1$ is sufficiently general, the intersection $S$ is irreducible by Bertini theorem. Then $S$ is an irreducible surface in the threefold $T_0$, dominant on both factors in the product $\PR ^2\times \PR ^2$.

Let
  $$
  Q_i=\nu (P_i)
  $$
be the images of the singular points $P_i$ in the threefold $T_0$. As meeting a point is a closed condition, the general section $T_1$ does not meet the points $Q_i$ in $\PR ^2\times \PR ^2$. Hence, the same property holds true for the surface $S$. 

Let 
  $$
  \mu =\nu |_V:V\to T_0
  $$
be the restriction of the morphism $\nu $ to $V$. Then $\mu $ is a finite surjective morphism from $V$ on to $T_0$. 

Let 
  $$
  F\subset \mu ^{-1}(S)
  $$ 
be an irreducible component of the scheme-theoretic preimage $\mu ^{-1}(S)$ of the surface $S$ under the morphism $\mu $. Then $F$ is an irreducible surface in $V$,
  $$
  F\cap V_{\sing }=\emptyset \; ,
  $$
and $F$ is dominant on both factors, i.e. $F$ is a fat surface in $X\times X$ which is away of the singularities of $V$.

Now, let $\phi $ be the generic point of $F$, and let
  $$
  K=\kappa (\phi )
  $$
be the field of rational functions on the surface $F$ over $k$. Consider the closed immersion
  $$
  F\to X\times X
  $$
and its composition
  $$
  F\to X
  $$
with the projection on to the first factor. This is a finite surjective morphism from $F$ to $X$.

The obvious commutative diagram 

\medskip

    $$
    \xymatrix{
    F\ar@/_/[dddr]_{} \ar@/^/[drrr]^-{\id _F}
    \ar@{.>}[dr]^-{} \\
    & X\times F \ar[dd] \ar[rr] & & F \ar[dd] \\ \\
    & X \ar[rr] & & \Spec (k)}
    $$
    
\medskip

\noindent gives us the graph embedding of $F$ in to $X\times F$ over $k$. Take the pullback of the top triangle to the generic point $\phi $ as shown here:

\medskip

    $$
    \xymatrix{
    \phi \ar[dd] \ar@/^/[drrr]^-{\id _{\phi }}
    \ar@{.>}[dr]^-{} \\
    & X\times \phi \ar[dd] \ar[rr] & & \phi \ar[dd] \\
    F\ar@/_/[dddr]_{} \ar@/^/[drrr]^-{\id _F}
    \ar@{.>}[dr]^-{} \\
    & X\times F \ar[dd] \ar[rr] & & F \ar[dd] \\ \\
    & X \ar[rr] & & \Spec (k)}
    $$
    
\medskip

\noindent Since 
  $$
  \phi =\Spec (K)\; ,
  $$
this pullback gives us a closed $K$-rational point 
  $$
  P_F:\Spec (K)\to X_K
  $$ 
on the surface $X_K$. 

Fix an embedding 
  $$
  K\subset \ud 
  $$
over $k$. Then $P_F$ can be also considered as a closed point in $X_{\ud }$. Since the surface $F$ is fat inside $X\times X$,
  $$
  \trdeg (P_F/k)=2\; .
  $$
  
\subsection{The Brevik-Nollet theorem}
\label{BrevikNolletThm}

Let now $k$ be the field of complex numbers $\CC $. Then, for a sufficiently big degree $d$, and for a sufficiently general section $V$, we have an isomorphism
   $$
   CH^1(V)=\langle \Delta _X\rangle \oplus I\; ,
   $$
where the first summand is the group generated by the class of the diagonal $\Delta _X$, and $I$ is the image of the restriction homomorphism 
   $$
   CH^1(W)\to CH^1(V)\; ,
   $$
where
  $$
  W=X\times X\; .
  $$
This is a particular case of the result of Brevik and Nollet, see Theorem 1.1 on page 108 in \cite{BrevikNollet}. 

Our surface $F$ is then rationally equivalent on $V$ to the linear combination 
  $$
  a\Delta _X+B_0 \; ,
  $$
where $B_0$ is the restriction of a divisor on $W$ to $V$. 

Assume that
  $$
  q=H^1(X,\bcO _X)=0
  $$
for the surface $X$. Then 
  $$
  CH^1(W)=CH^1(X)\oplus CH^1(X)\; ,
  $$
see Exercise 12.6 on page 292 in \cite{Hartshorne}. It follows immediately, that $B_0$ is the restriction of a linear combination of divisors of type $D\times X$ and $X\times D$ from $W=X\times X$ to the threefold $V$. As such, $B_0$ is a balanced divisor on $V$, i.e. a $2$-cycle balanced on $X\times X$ and supported on $V$. 

Moreover, since any curve on $X$ can be moved away from a finite collection of points in its rational class on $X$, by the Chow moving lemma, we can also assume that the support of $B_0$ does not meet the singular points $P_i$. 

Notice also that, if $k$ is any characteristic $0$ algebraically closed field of definition, we still can apply the Brevik-Nollet theorem by embedding the finitely generated minimal field of definition in to $\CC $ (the Lefschetz principle).

\medskip

There are two different cases here.

\medskip

\noindent {\it Case 1: $a=0$}

\medskip

Then $F$ is balanced modulo rational equivalence on $V$, and so on $X\times X$. It means that the closed point $P_F$ is rationally equivalent to a point of transcendence degree $<2$ on $X_{\ud }$. Then $\Delta _X$ is balanced on $X\times X$, and the motive $M(X)$ is a submotive in the direct sum of motives of curves by Theorem \ref{maintrdeg}. Then, of course, $M(X)$ is finite-dimensional by Proposition \ref{Kimuraab}.

\bigskip

\noindent {\it Case 2: $a\neq 0$}

\medskip

In this case the diagonal $\Delta _X$ is rationally equivalent, with coefficients in $\QQ $, to the divisor
  $$
  \frac{1}{a}\, (F-B_0)
  $$
on $V$, where both $F$ and $B_0$ do not meet the singular locus $V_{\sing }$. 

\subsection{The application of Weak Lefschetz}
\label{WL}

Now assume that the group $H^*(X,\QQ _l)$ is algebraic. This is our fundamental assumption, the assumption of Bloch's conjecture. Recall that it is equivalent to
  $$
  H^1(X,\QQ _l)=H^3(X,\QQ _l)=0
  $$
and
  $$
  H^2_{\tr }(X,\QQ _l)=0\ .
  $$
Since $k=\CC $, vanishing of odd cohomology implies $q=0$. 

As we explained in Section \ref{motivic-II}, BC is equivalent to proving that the diagonal $\Delta _X$ is a balanced correspondence on $X\times X$. 

In Case 1, when the coefficient $a$ is $0$, the diagonal is balanced, the motive $M(X)$ is a submotive in a direct sum of motives of curves, and we are done by Proposition \ref{Kimuraab}. Suppose we are in Case 2, that is, $a\neq 0$. 

Let $V_{\sm }$ be the smooth locus of the threefold $V$, and let $Z$ be an integral closed subscheme of codimension $1$ in $V$, such that the support of $Z$ does not meet $V_{\sing }$. Then we can define the cohomological class $\cl ^1_V(Z)$ of $Z$ in $H^2(V,\QQ _l)$ using the excision isomorphism between local cohomology groups  $H^2_Z(V_{\sm },\QQ _l)$ and $H^2_Z(V,\QQ _l)$ and the Thom isomorphism, see page 605 in \cite{Srinivas}. Extending by linearity, we have a legal cohomology class $\cl ^1_V(Z)$ in $H^2(V,\QQ _l)$ for each divisor $Z$ whose support does not meet the singular locus $V_{\sing }$.

In particular, since both $F$ and $B$ do not meet the singular locus $V_{\sing }$, we have the cohomology classes 
  $$
  \cl ^1_V(F)\qqand \cl ^1_V(B_0)
  $$
in the group $H^2(V,\QQ _l)$.

The restriction homomorphism
   $$
   H^2(X\times X,\QQ _l)\to H^2(V,\QQ _l)
   $$
 is an isomorphism by Weak Lefschetz\footnote{for the \'etale version of the Weak Lefschetz Theorem see, for example, Corollary 9.4 on page 106 in \cite{FreitagKiehl} or \cite{SGA7-2}, Expos\'e XVII and XVIII}. Since odd chomology groups vanish for $X$, by the K\"unneth formula, we obtain that
  $$
  H^2(X\times X,\QQ _l)=H^2(X,\QQ _l)\oplus H^2(X,\QQ _l)\; .
  $$
  
Let $b_2$ be the second Betti number of the surface $X$, and let
  $$
  D_1,\ldots ,D_{b_2}
  $$
be divisors on $X$ generating $H^2(X,\QQ _l)$. Then the group
  $$
  H^2(X\times X,\QQ _l)
  $$
is generated by the cohomological classes of the divisors $D_i\times X$ and $X\times D_j$, for all $i$ and $j$.
  
Then we obtain that the class $\cl ^1_V(F)$ is represented by a linear combination
  $$
  B_1=\sum _im_i(D_i\times X)|_V+
  \sum _jn_j(X\times D_j)|_V
  $$
of the restrictions of such divisors to the threefold $V$. Clearly, these restrictions are balanced, and therefore $B_1$ is balanced on $V$ as a threefold in $X\times X$. 

\subsection{Matsusaka's coincidence}
\label{matsusaka}

Let
  $$
  \tilde V\to V
  $$
be the resolution of singularities $P_i$ on then threefold $V$. This resolution can be constructed as follows. 

Write again
  $$
  W=X\times X\; ,
  $$
and consider the composition
  $$
  \tilde W\to W
  $$
of the consequent blowups of $W$ at the points $P_i$. Since $W$ is smooth, exceptional fibres are the projective spaces $\PR ^3$. Then $\tilde V$ can be viewed as the strict transform of $V$ in $\tilde W$, and the exceptional divisor at $P_i$ is a smooth quadric 
  $$
  E_i\subset \PR ^3\; .
  $$

Let 
  $$
  E=\cup E_i
  $$
be the union of these quadrics in $\tilde V$, and let
  $$
  U=\tilde V\smallsetminus E=V\smallsetminus V_{\sing }
  $$
be the complement to the exceptional locus in $\tilde V$, or to the singular locus in $V$.

Since $F$ does not meet the points $P_i$, we have that
  $$
  F\subset U\; .
  $$
The variety $U$ is smooth, and the restriction of $B_1$ on $U$ is cohomologically equivalent to $F$, that is, 
  $$
  \cl ^1_U(F)=\cl ^1_U(B_1|_U)\; .
  $$ 
 
Let
  $$
  \tilde B_1
  $$
be the strict transform of the balanced divisor $B_1$ on $V$ in the smooth threefold $\tilde V$. Then we also have that $F$ is cohomologically equivalent to $\tilde B_1|_U$ on $U$. And then the long exact sequence 
  $$
  \ldots \to H^2_E(\tilde V,\QQ _l)\to H^2(\tilde V,\QQ _l)\to 
  H^2(U,\QQ _l)\to \ldots 
  $$
gives us that
  $$
  F-\tilde B_1
  $$
is cohomologically equivalent to a linear combination
  $$
  \sum _im_iE_i
  $$
on $\tilde V$.

Now recall the Matsusaka theorem, which says that cohomological equivalence relation coincides with the algebraic one, for divisors with coefficients in $\QQ $ on smooth projective varieties over an algebraically closed field of any characteristic, see Theorem 4 on page 65 in \cite{Matsusaka}. 

As $\tilde V$ is a smooth projective variety over an algebraically closed field, we can apply the Matsusaka theorem and obtain that the difference $F-\tilde B_1$ is algebraically equivalent, with rational coefficients, to the linear combination $\sum _im_iE_i$ on $\tilde V$.
 
And since proper push-forwards preserve algebraic equivalence, we see that $F-\tilde B_1$ is algebraically equivalent to $0$, that is $F$ is algebraically equivalent to $B_1$ on $V$, and hence on $X\times X$, with coefficients in $\QQ $. 

Then the diagonal
  $$
  \Delta =\Delta _X
  $$ 
is algebraically equivalent to the balanced $2$-cycle
  $$
  \gB =\frac{1}{a}\, (B_1-B_0)
  $$
on $X\times X$ with coefficients in $\QQ $. 
 
\subsection{Coda}

Then 
  $$
  \gA =\Delta -\gB 
  $$ 
is a smash-nilpotent correspondence on $X$ by Theorem \ref{VNT1}. 

By Corollary \ref{VNT2}, $\gA $ is also a nilpotent element in the associative algebra of degree $0$ correspondences from $X$ to $X$, i.e.
  $$
  \gA ^n=(\Delta -\gB )^n=0
  $$
in  
  $$
  CH^2(X\times X,\QQ )
  $$ 
for some positive integer $n$. 

As $\Delta $ is the identity element for multiplication, Newton's binomial formula yields 
  $$
  \gA ^n=\sum _{i=0}^n{n\choose i}(-1)^{n-i}\gB ^{n-i}=0\; ,
  $$
whence 
  $$
  \Delta =\sum _{i=0}^{n-1}{n\choose i}(-1)^{n-i+1}\gB ^{n-i}
  $$
is a balanced correspondence in $CH^2(X\times X,\QQ )$. 

Now, as $\gB $ is balanced, so is the diagonal $\Delta $. And since $\Delta =\Delta _X$ is balanced, the motive $M(X)$ is a submotive in the direct sum of motives of curves by Theorem \ref{maintrdeg}. Hence $M(X)$ is in $\CHM (k,\QQ )^{\ab }$ and, therefore, finite-dimensional by Proposition \ref{Kimuraab}. Then Conjecture \ref{BCgeneralized} is true by Corollary \ref{Kimurakey2}. 

Theorem A is proved.

Notice that we used Brevik-Nollet's theorem, which is, so far, proven over $\CC $ only (see, however, \cite{LenaJi}). However, the Lefschetz principle and Lemma \ref{scalars} guarantee that our proof extends to any algebraically closed field $k$ of characteristic $0$. 

\section{When $H^2$ is not algebraic}
\label{notalg}

If $k=\CC $ and $H^2(X,\QQ _l)$ is not algebraic, that is $p_g\neq 0$, then the irregularity $q=0$ and Case 1 are not possible simultaneously by Proposition 2.4 in \cite{Barbieri-Viale} (see also Mumford's result in \cite{Mumford}).

Another temptation would be to use Betti cohomology and Hodge decomposition to argue as follows. 

Assuming $q=0$ and $p_g>0$, let
  $$
  \phi =\cl ^1_V(F)\; ,
  $$
and let 
  $$
  \psi \in H^2(X\times X,\QQ )
  $$
be a class, such that
  $$
  \Res (\psi )=\phi 
  $$
under the restriction homomorphism 
  $$
  \Res :H^2(X\times X,\QQ )\to H^2(V,\QQ )\; ,
  $$
all in terms of Betti cohomology with coefficient in $\QQ $ now. 

Tensoring with $\CC $, we obtain the complexified cohomology classes 
  $$
  \psi _{\CC }\in H^2(X\times X,\CC )
  $$
and 
  $$
  \phi _{\CC }\in H^2(V,\CC )
  $$
respectively, and
  $$
  \phi _{\CC }=\Res (\psi _{\CC })\; ,
  $$
under the restriction on cohomology with coefficients in $\CC $. 

Let
  $$
  \psi _{\CC }=
  \psi ^{2,0}_{\CC }+\psi ^{1,1}_{\CC }+\psi ^{0,2}_{\CC }
  $$
be the sum induced by the Hodge decomposition
  $$
  H^2(X\times X,\CC )=H^{2,0}(X\times X)\oplus H^{1,1}(X\times X)\oplus H^{0,2}(X\times X)\; .
  $$
Since $\Res $ is a morphism of Hodge structures, we obtain that
  $$
  \Res (\psi ^{1,1}_{\CC })=\phi _{\CC }
  $$
and
  $$
  \Res (\psi ^{2,0}_{\CC })=\Res (\psi ^{0,2}_{\CC })=0\; ,
  $$
where $\Res $ is the complexified restriction 
  $$
  \Res :H^2(X\times X,\CC )\to H^2(V,\CC )\; .
  $$

The temptation would be now to think that $\psi ^{1,1}_{\CC }$ is in the intersection
  $$
  H^{1,1}(X\times X)\cap H^2(X\times X,\QQ )\; ,
  $$
and so algebraic by $(1,1)$-Lefschetz theorem. 

However, if that could be true, then $\phi $ would be balanced, because
  $$
  CH^1(X\times X,\QQ )=CH^1(X,\QQ )\oplus CH^1(X,\QQ )
  $$
thanks to vanishing $H^1(X,\QQ )=0$. This, again, contradicts to Proposition 2.6 in \cite{Barbieri-Viale} and Mumford's result in \cite{Mumford}. 

This means that the class $\psi ^{1,1}_{\CC }$, being a summand in a linear transformation of the complexified class $\psi _{\CC }$ inside the complex vector space $H^{1,1}(X\times X)$, is not itself a complexification of a cycle class in $H^2(X\times X,\QQ )$.

The study of the class $\psi ^{1,1}_{\CC }$ may well be important in approaching the problem of motivic finite-dimensionality of $K3$-surfaces over $\CC $.


\section{Proof of Theorem B}

\subsection{Over $\bar \FF _p$}

To prove Bloch's conjecture in characteristic $p>0$, we first repeat all the same arguments until the point when we need to apply the Brevik-Nollet theorem, which is not yet proven in positive characteristic. 

In particular, we construct the fat surface $F$ on $V$, away from $V_{\sing }$, by the method presented in Section \ref{fat}. Since $F\cap V_{\sing }=\emptyset$, we have an honest cohomological class 
  $$
  \cl ^1_V(F)
  $$ 
of $F$ in $H^2(V,\QQ _l)$ by excision. By Weak Lefschetz, $\cl ^1_V(F)$ is coming from 
  $$
  H^2(X\times X,\QQ _l)
  $$ 
 via restriction from $X\times X$ to the threefold $V$. 
 
 Since $H^*(X,\QQ _l)$ is algebraic, any element in $H^2(X\times X,\QQ _l)$ is the class of a linear combination of divisors of type $D\times X$ or $X\times D$, where $D$ is a divisor on $X$. Therefore, $F$ is cohomologically equivalent to a balanced divisor $B$ not meeting singular points on $V$. 

As above, let $\tilde V$ be the resolution of singularities on $V$. Since $F$ and $B$ do not meet $V_{\sing }$, we can consider them as divisors on $\tilde V$. By Matsusaka's theorem, $F$ is algebraically equivalent to $B$ on $\tilde V$. 

Now let 
  $$
  k=\bar \FF _p\; .
  $$
Then, of course, the Albanese kernel vanishes over $k$, see Corollary 3.3 on page 324 in \cite{Raskind} and Proposition 4 on page 788 in \cite{CTSS}. 

But we can prove more. 

\subsection{Over $\overline {\FF _p((t))}$}

Let
  $$
  \ud =\overline {\FF _p((t))}
  $$
be the closure of Laurent series with coefficients in $\FF _p$. This is an uncountable universal domain in the sense of Andr\'e Weil, which may be also considered as the field of generalised Puiseux series over $\FF _p$, see \cite{Kedlaya}. Let us prove that the Albanese kernel vanishes also for $X_{\ud }$. 

Indeed, consider the Picard variety 
  $$
  A=\Pic ^0_{\tilde V/k}
  $$
of the smooth projective threefold $\tilde V$ over $k$, and apply the old argument as in the proof of Proposition 4 on page 788 in \cite{CTSS}. 

Namely, since the divisor $F-B$ is algebraically equivalent to $0$ on $\tilde V$, it gives the point $P$ in $A(\bar \FF _p)$. As the latter group is torsion, it follows that $F$ is rationally equivalent to $B$ with coefficients in $\QQ $ on $\tilde V$. 

Now again, let $K$ be the field of rational functions on $F$, and let $P_F$ be the closed point in $X_K$ induced by the composition of the embedding of $F$ in to $X\times X$ and the projection on to the first factor. Since $F$ is fat, 
  $$
  \trdeg (P_F/k)=2\; ,
  $$
and since $B$ is balanced, 
  $$
  \trdeg ([P]/k)<2\; .
  $$ 
By Theorem \ref{maintrdeg}, the diagonal $\Delta _X$ is balanced, the motive $M(X)$ is finite-dimensional, and 
  $$
  T(X_{\ud })=0
  $$ 
by Corollary \ref{Kimurakey2}.

\section{After proof}

\subsection{Finite dimension of $M(X)$}

Let again $k$ be an algebraically closed field of $\cha (k)\neq 2$, and let $X$ be a smooth projective surface with algebraic $H^*(X,\QQ _l)$ over $k$. Let $M(X)$ be the Chow motive of the surface $X$, and let 
  $$
  M(X)=\uno \oplus M^2(X)\oplus \Le ^2
  $$
be the decomposition of $M(X)$ in the sense of Murre, see \cite{Murre}. Here $\uno $ is the unit motive and $\Le ^2$ is the square of the Lefschetz one. 

Notice that the Picard motive $M^1(X)$ and the Albanese motive $M^3(X)$ are both of abelian type, and hence finite-dimensional by Proposition \ref{Kimuraab}. Since, moreover, 
  $$
  H^1(X,\QQ _l)=H^3(X,\QQ _l)=0\; ,
  $$
we obtain that
  $$
  M^1(X)=M^3(X)=0\; ,
  $$
because finite-dimensional motives cannot be phantoms, see Corollary 7.3 on page 192 in \cite{Kimura}. 

Now, the motive $M(X)$ is finite-dimensional by the proven Conjecture \ref{BCgeneralized}, Theorem \ref{maintrdeg} and Proposition \ref{Kimuraab}. And since $M^1=M^3=0$, the motive $M(X)$ is evenly finite-dimensional. That is,
  $$
  \wedge ^{b+1}M(X)=0\; ,
  $$
where 
  $$
  b=\dim (H^*(X))
  $$ 
is the total Betti number of the surface $X$. 

\subsection{The motive $M^2_{\tr }$}

Let also
  $$
  M^2_{\tr }(X)
  $$
be the transcendental motive of $X$, as defined in \cite{KahnMurrePedrini}. It is a submotive in $M^2(X)$, whose second cohomology group vanishes,
  $$
  H^2(M^2_{\tr }(X))=H^2_{\tr }(X,\QQ _l)=0\; ,
  $$
if $H^2(X,\QQ _l)$ is algebraic. Here $H^2_{\tr }(X,\QQ _l)$ is the second transcendental cohomology group in terms of Bloch, see Section \ref{nonmotivic}.

Since the motive $M(X)$ is finite-dimensional, so is its submotive $M^2_{\tr }(X)$. Then
  $$
  M^2_{\tr }(X)=0
  $$
by the same Corollary 7.3 in \cite{Kimura}.  

\bigskip

\bigskip

\begin{small}

\end{small}

\bigskip

\bigskip

\begin{small}

{\sc Department of Mathematical Sciences, University of Liverpool, Peach Street, Liverpool L69 7ZL, England, UK}

\end{small}

\medskip

\begin{footnotesize}

{\it E-mail address}: {\tt vladimir.guletskii@liverpool.ac.uk}

\end{footnotesize}


\begin{thebibliography}{9999999}







\bibitem[SGA$4\frac{1}{2}$]{SGA4.5}
Cohomologie \'etale. S\'eminaire de G\'eom\'etrie Alg\'ebrique du Bois Marie (SGA $4\frac{1}{2}$) dirig\'e par P. Deligne avec la collaboration de J.F. Boutot, A. Grothendieck, L. Illusie et J.L. Verdier. Lecture Notes in Mathematics. Vol. 569 (1977)


\bibitem[SGA7-2]{SGA7-2}
Groupes de monodromie en g\'eom\'etrie alg\'ebrique. S\'eminaire de G\'eom\'etrie Alg\'ebrique du Bois Marie 1967-69 (SGA 7-2) par P. Deligne, N. Katz. Lecture Notes in Mathematics 340 (1973)

\bibitem[Stacks]{StacksProject}
The Stacks Project. \url{https://stacks.math.columbia.edu}




\bibitem{AK}
Y. Andre, B. Kahn (with appendix by P. O'Sullivan). Nilpotence,
radicaux et structure monoidales. Rend. Sem. Mat. Univ. Padova 108 (2002) 107 - 291



\bibitem{Barlow1}
R. Barlow. A simply connected surface of general type with $p_g=0$. Inventiones mathematicae 79 (1985) 293 - 301

\bibitem{Barlow2}
R. Barlow. Rational equivalence of zero cycles for some more surfaces with $p_g=0$. Inventiones mathematicae 79 (1985) 303 - 308


\bibitem{Barbieri-Viale}
L. Barbieri-Viale. Balanced varieties. Algebraic $K$-theory and its applications (Trieste, 1997), 298 - 312, World Sci. Publ., River Edge, NJ, 1999




\bibitem{BlochThesis}
S. Bloch. Algebraic Cohomology Classes on Algebraic Varieties. PhD thesis 1970

\bibitem{BCinitial}
S. Bloch. $K_2$ of artinian $\QQ $-algebras, with application to algebraic cycles. Communications in Algebra. Volume 3, Issue 5 (1975) 405 - 428 

\bibitem{BlochLectures}
S. Bloch. Lectures on Algebraic Cycles. Duke University Math. Series IV. Durham, NC. Duke University 1980

\bibitem{BKL}
S. Bloch, A. Kas, D. Lieberman. Zero-cycles on surfaces with $p_g=0$. Compositio Math. Tome 33 (1976) 135 - 145

\bibitem{BlochSrinivas}
S. Bloch, V. Srinivas. Remarks on correspondences and algebraic cycles. Amer. J. Math. 105 (1983) 1235 - 1253








\bibitem{BrevikNollet}
J. Brevik. S. Nollet. Grothendieck-Lefschetz theorem with base locus. Israel Journal of Mathematics 212 (2016) 107 - 122




\bibitem{CTSS}
J.-L. Colliot-Th\'el\`ene, J.-J. Sansuc, C. Soul\'e. Torsion dans le groupe de Chow en codimension deux. Duke Mathematical Journal. Vol 50 (1983) 763 - 801 




\bibitem{DeligneMilne}
P. Deligne and J. Milne. Tannakian categories. In Hodge Cycles and Shimura Varieties. Lecture Notes in Math. 900 (1982) 101 - 208



\bibitem{FreitagKiehl}
E. Freitag, R. Kiehl. Etale Cohomology and the Weil Conjecture. Ergebnisse der Mathematik und ihrer Grenzgebiete. 3 Folge. A Series of Modern Surveys in Mathematics. Volume 13. Springer-Verlag 1988


\bibitem{Fulton}
W. Fulton. Intersection theory. Ergebnisse der Mathematik und ihrer Grenzgebiete. 3 Folge. Band 2. Springer-Verlag 1984




\bibitem{trdegzerocycles}
S. Gorchinsky, V. Guletski\u \i . Transcendence degree of zero-cycles and the structure of Chow motives. Central European J. Math. Vol. 10, No. 2 (2012) 559 - 568


\bibitem{GP2}
V. Guletski\u {\i }, C. Pedrini. Finite-dimensional motives and the conjectures of Beilinson and Murre. $K$-Theory. Vol. 30, No. 3 (2003) 243 - 263

\bibitem{Hartshorne}
R. Hartshorne. Algebraic Geometry. Springer 1977






\bibitem{InoseMizukami}
H. Inose, M. Mizukami. Rational equivalence of zero-cycles on some surfaces with $p_g=0$. Math. Ann. 244 (1979) 205 - 217







\bibitem{LenaJi}
L. Ji. The Noether-Lefschetz theorem in arbitrary characteristic. Journal of Alg. Geometry. Vol. 33, No. 3 (2024) 567 - 600

\bibitem{KahnMurrePedrini}
B. Kahn, J. Murre, C. Pedrini. On the transcendental part of the motive of a surface. Algebraic cycles and motives. Vol. 2. London Math. Soc. Lecture Note Ser. Vol. 344 (2007) 143 – 202. Cambridge University Press

\bibitem{Kedlaya}
K. Kedlaya. The algebraic closure of the power series field in positive characteristic. Proc. Amer. Math. Soc. 129 (12) (2001) 3461 - 3470

\bibitem{Kimura}
S.-I. Kimura. Chow groups are finite dimensional, in some sense. Math. Ann. Vol. 331 (2005) no. 1, 173 - 201


\bibitem{Kleiman}
S. Kleiman. Algebraic cycles and the Weil conjectures. Dix espos\' es sur la cohomologie des sch\' emas. Amsterdam, North-Holland (1968) 359 - 386






\bibitem{LiedtkeGodeaux}
C. Liedtke. Non-classical Godeaux surfaces. Math. Ann. Vol. 343,  No. 3 (2009) 623 - 637

\bibitem{Liedtke}
C. Liedtke. Algebraic Surfaces in Positive Characteristic. In Birational Geometry, Rational Curves, and Arithmetic. Simons Symposia. Springer (2013) 229 - 292






\bibitem{Matsusaka}
T. Matsusaka. The Criteria for Algebraic Equivalence and the Torsion Group. American J. of Math. Vol. 79, No. 1 (1957) 53 - 66


\bibitem{MilneEC}
J. Milne. \'Etale Cohomology. Princeton University Press 1980


\bibitem{MilneRoitmanThm}
J. Milne. Zero cycles on algebraic varieties in nonzero characteristic: Rojtman's theorem. Compositio Mathematica. Vol. 47, No. 3 (1982) 271 - 287


\bibitem{Mumford}
D. Mumford. Rational equivalence of $0$-cycles on surfaces. J. Math. Kyoto Univ. 9 (1969) 195 - 204


\bibitem{Murre}
J. Murre. On the motive of an algebraic surface. J. f\"{u}r die reine und angew. Math. 409 (1990) 190 - 204







\bibitem{Raskind}
W. Raskind. Algebraic K-theory, \'etale cohomology and torsion algebraic cycles. Contemporary Mathematics. Vol. 83 (1989) 311 - 341

\bibitem{Roitman}
A. Roitman. The torsion of the group of $0$-cycles modulo rational equivalence. Annals of Mathematics. Second Series 111 (3) (1980) 553 - 569





\bibitem{Scholl}
A. Scholl. Classical motives. In "Motives", Proc. Symp. Pure Math. Vol.55, Part 1 (1994) 163 - 187








\bibitem{Srinivas}
V. Srinivas. Algebraic Cycles on Singular Varieties. Proceedings of the International Congress of Mathematicians in Hyderabad, India (2010) 603 - 623





\bibitem{VoevodskyNilp}
V. Voevodsky. A nilpotence theorem for cycles algebraically equivalent to zero. Int. Math. Res. Notices. Vol. 1995, Issue 4 (1995) 187 - 199


\bibitem{Sur les 0-cycles}
C. Voisin. Sur les z\'ero-cycles de certaine hypersurfaces munies d'un automorphisme. Ann. Scuola Norm. Sup. Pisa Ck. Si. (4) 19 (1992) 473 - 492

\bibitem{VoisinVariations}
C. Voisin. Variations de structure de Hodge et z\'ero-cycles sur les surfaces g\'en\'erales. Math. Ann. 299 (1994) 77 - 103

\bibitem{VoiN}
C. Voisin. Remarks on zero-cycles of self-products of varieties. In: Moduli of vector bundles. Lecture Notes in Pure and Appl. Math. 179. Dekker, New York (1996) 265 - 285

\bibitem{VoisinCataneseBarlowSurfaces}
C. Voisin. Bloch's conjecture for Catanese and Barlow surfaces. J. Differential Geom. Volume 97, Number 1 (2014) 149 - 175





\end{thebibliography}
\end{document}